\documentclass{comnet}
\usepackage[utf8]{inputenc}
\usepackage{amsmath}
\usepackage{mathtools}
\usepackage{support-caption}
\usepackage{subcaption}
\usepackage{wrapfig}
\usepackage{float}
\usepackage{tcolorbox}

\DeclarePairedDelimiter\fl{\lfloor}{\rfloor} 
\DeclarePairedDelimiter\cb{\{}{\}} 

\begin{document}

\title{Patterns of Primes and Composites from Divisibility Network of Natural Numbers}

\shorttitle{Divisibility Network of Natural Numbers} 
\shortauthorlist{Abiya Rajans, G. Ambika} 

\author{
\name{Abiya Rajans} and \name{G. Ambika{$^*$}}
\address{Indian Institute of Science Education and Research (IISER) Tirupati, Tirupati - 517 507, India\email{$^*$Corresponding author: g.ambika@iisertirupati.ac.in}}}

\maketitle

\begin{abstract}
{We present the pattern underlying some of the properties of natural numbers, using the framework of complex networks. The network used is a divisibility network in which each node has a fixed identity as one of the natural numbers and the connections among the nodes are made based on the divisibility pattern among the numbers.  We derive analytical expressions for the centrality measures of this network in terms of the floor function and the divisor functions. We validate these measures with the help of standard methods which make use of the adjacency matrix of the network. Thus how the measures of the network relate to patterns in the behaviour of primes and composite numbers becomes apparent from our study.} 
{Complex networks, divisibility network, primes and composites,  stretching similarity, floor function, divisor function}
\\
\end{abstract}

\section{Introduction}

In the world of numbers, primes and composites form two non-overlapping infinite
sets. Composites can be factorised into primes, while primes remain elusive in many
aspects. Prime numbers are found to have interesting connections underlying widely
different phenomena and processes. They repeatedly pop up in many unexpected
ways in Riemann $\zeta$ function, self-organised criticality, quantum computation,
quantum cosmology, signal processing, fractal geometry and public-key cryptography
algorithms. \citep{Hardy1975,Luque2008, Shor1994, Sanchis2012,Dragovich2010, Cattani2010, Riesel1994, Cusick2004}. The positions of prime numbers appear to be scattered throughout natural numbers in a non-homogeneous
random fashion and the most famous theorem that relates this is the Prime number
theorem which states that the number of primes up to $N$, approaches $\frac{N}{\log{N}}$, in the
limit of $N$ $\rightarrow \infty$ \citep{Hardy1975}.

Attempts to understand the architecture of natural numbers have led to much
exciting research, involving methods from various fields ranging from number
theory to graph theory. The framework of complex networks has been used
recently to understand the pattern of natural numbers in three different ways of
construction. One of the studies reported relies on constructing a weighted bipartite
of composite and prime numbers, with connections between them decided by the
prime factorisation \citep{Garcia2014}. Similarly, another work studies a network of composite numbers where the link between two nodes exists if they share a common divisor \citep{Lewis2008}. The congruence
relations among numbers are explored using a multiplex network, and every layer is
reported as a sparse and heterogeneous sub-network with a scale-free topology \citep{XiaoYong2016}.
A more general approach is to consider all the numbers as nodes, and if one number
divides another, their corresponding nodes are connected by an undirected link. This
results in a deterministic but complex network called natural number
network, which has been reported to be scale-free with an interesting
symmetry property called stretching similarity \citep{Shekatkar2015}. Following this, a recent work has
applied this to the specific case of divisibility pattern within the elements of Pascal
triangles \citep{Solares2020}. In all the above cases, the topology of the resulting networks is shown
to be of scale-free type. 

We note that approaches based on the framework of complex networks are ideal to
understand and visualise the hidden structure in the sequence of primes and their
relations with other natural numbers. Since divisibility pattern is the one aspect that
distinguishes the primes from other numbers, we take the measures of natural number network constructed using divisibility patterns to discern the intricacies in the architecture of primes and
composites in a natural way. In this work, we present how the measures of this   network can
be related to the properties of natural numbers, especially prime numbers. Here, the
pattern of connectivity reflects the divisibility pattern among the numbers and hence
can project primes very effectively. Moreover, as an advancement of earlier work,
we derive analytic expressions for the network measures like degree, clustering
coefficient, centrality and link density in terms of floor functions and divisor functions
and verify them by direct numerical simulations. These expressions help to
understand the trends observed in the characterising measures of the network as $N$
increases. As a specific case, we could explain the stretching similarity reported in
previous work \cite{Shekatkar2015} for local clustering coefficients of the nodes in this network as arising due to the trends in the divisor functions of natural numbers. We also bring out the specific
trends shown by the measures corresponding to prime numbers.

\section{Divisibility Network of Natural Numbers}

We define the divisibility network of natural numbers of size $N$ as a graph $G_N$ with vertices, $V(G_N) = \cb*{1,\dots,N}$ and edges connected according to the adjacency matrix, $A = [A_{ij}]$ where

$A_{ij}$ = 
$\begin{cases}
1, & \text{if $i \ne j$ and either $i$ divides $j$ or $j$ divides $i$} 
\\0, & \text{otherwise} 
\end{cases}$

\begin{figure}[hbtp]
\centering
\begin{subfigure}{0.49\textwidth}
\includegraphics[width=1\textwidth]{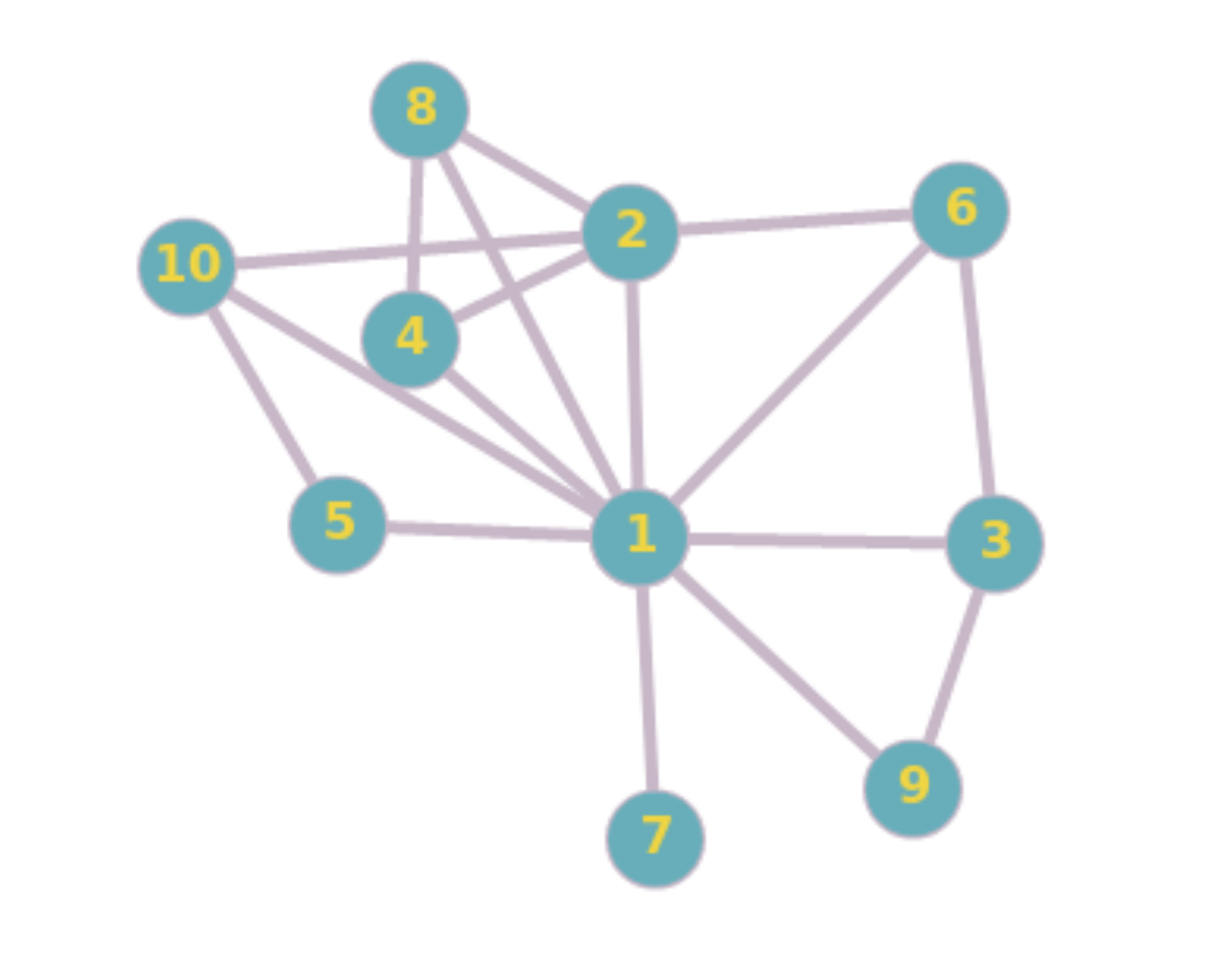}
\caption{  }\label{fig:1a}
\end{subfigure}
\vspace*{\fill} 
\begin{subfigure}{0.49\textwidth}
\includegraphics[width=1\textwidth]{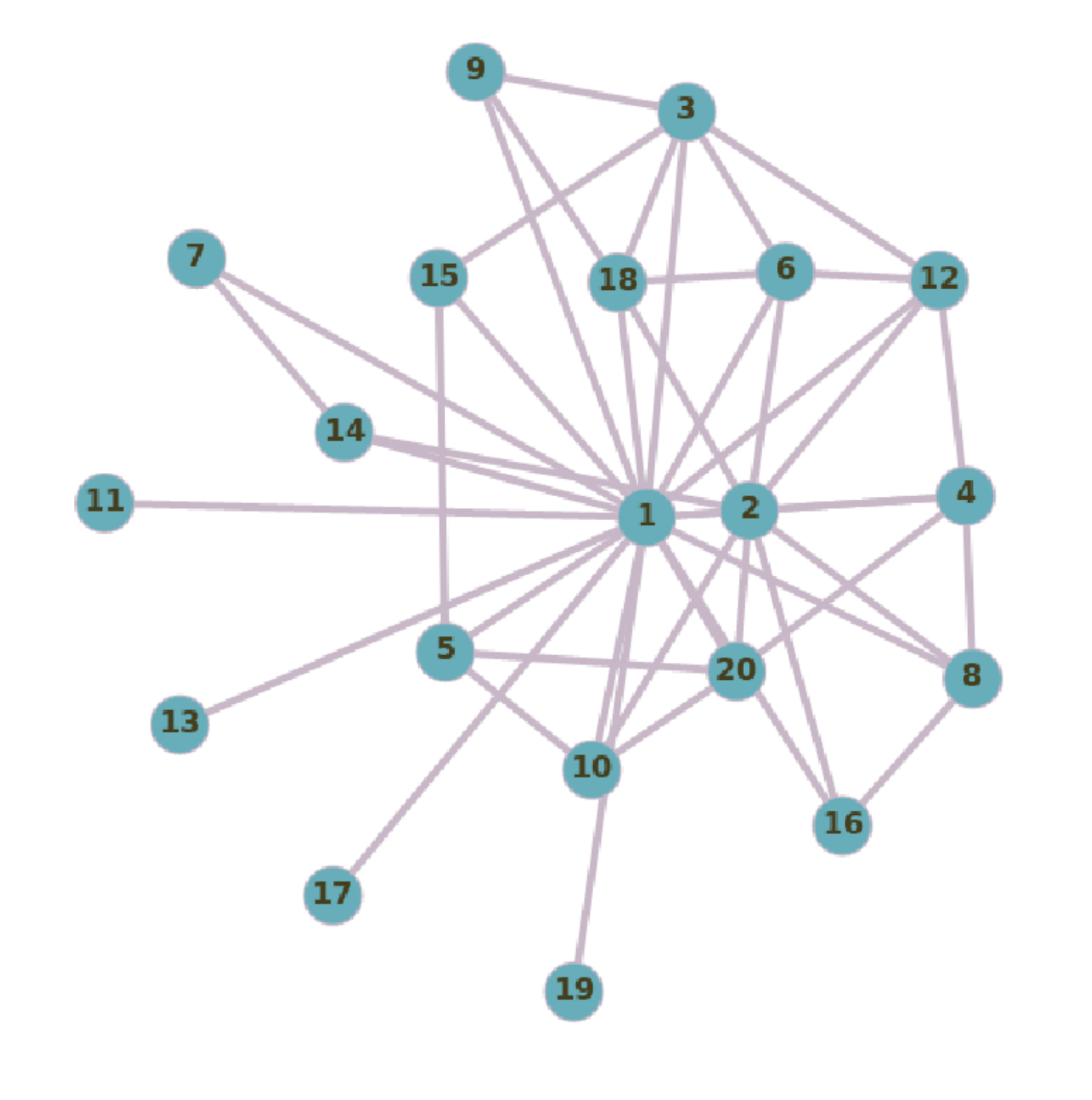}
\caption{ }\label{fig:1b}
\end{subfigure}
\caption{ Network of natural numbers $G_N$ of size $N$ with connectivity based on the divisibility patterns among the natural numbers. (a) $N$ = 10 and (b) $N$ = 20. }\label{fig:1}
\end{figure}

By construction, $G_N$ is a deterministic network whose nodes are fixed as natural numbers in order from 1 to $N$ and therefore is different from other networks like random or scale free networks. 
We present in figures \ref{fig:1a}, \ref{fig:1b} two such networks $G_N$ as illustrations, for $N$ = 10 and 20. 

We now derive expressions for some of the local characteristics of $G_N$ like degree of the nodes and related global measures and study their trends for prime and composite numbers.  

\section{Degrees of primes and composites on $G_N$ }

The connectivity of the network $G_N$ maps the pattern of divisibility among natural numbers. To characterise them, we start with the most basic local property of nodes in the network, namely the degree of a node. It is defined as the total number of links of the node \citep{Newman2010}. We study the degrees of all the nodes and their variation with the network size $N$. We find that we can express the degree of any node in terms of floor function and divisor function of that natural number. In particular, we first consider nodes that correspond to prime numbers among the first $N$ numbers. This will bring out patterns that are specific to primes. The expressions for the degrees derived below will help to compute the degree for any node on the network. Starting from the definition of degree, we can also numerically compute them from the elements of the adjacency matrix. We present both calculations to validate the expressions derived from the properties of numbers.

By the definition of degree as the number of links of a node, we get the degree of a node from the adjacency matrix of the constructed network as 

\begin{align}\label{3.1}
k_p = \sum_{j=1}^N A_{pj}
\end{align}

On the other hand, since $G_N$ is a deterministic network, the degree can be expressed as a function of the nodes. The degree of a prime node is the sum of its number of multiples and its only divisor 1. In a network of size $N$, $M(p) = \fl*{\frac{N}{p}}$ is the number of multiples of $p$ which are less than $N$, including $p$ itself. Since there are no self-loops in $G_N$, the multiples of $p$ less than $N$ excluding $p$ are $2.p, \  3.p, \  \dots, \  M(p).p$. Hence the node p is connected to $\fl*{\frac{N}{p}}-1$ multiples. Here $\fl*{\frac{N}{p}}$ known as the floor function, is the greatest integer value of $\frac{N}{p}$ \citep{Hardy1975}. For example, if $N$ = 5 and p = 2, then $\fl*{\frac{5}{2}} = \fl*{2.5}$ = 2. Now including the one divisor of $p$, the degree of p is $\fl*{\frac{N}{p}} - 1 + 1 = \fl*{\frac{N}{p}}$. 
Thus we see a trend in the degree of primes $k_{p}$ as a function of primes $p$ for a given network of size $N$, as 
\begin{align}\label{3.2}
   k_{p} = \fl*{\frac{N}{p}}
\end{align}

\begin{figure}[hbtp]
\centering
\includegraphics[width=1\textwidth]{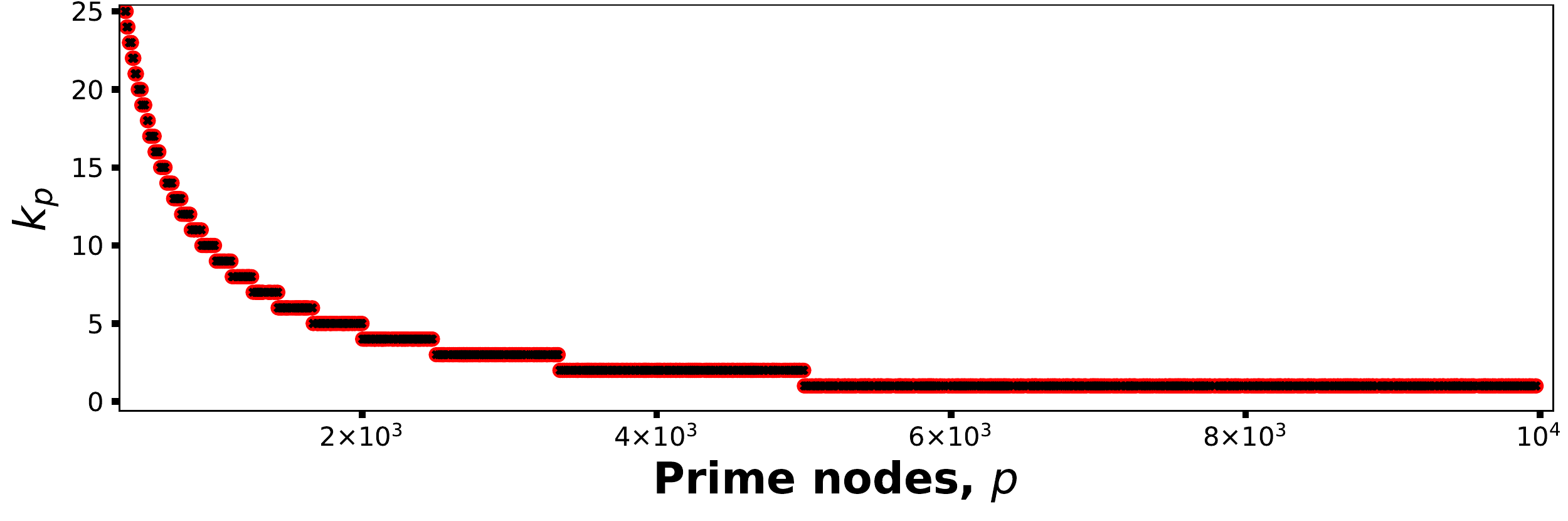}
\caption{(a) Degrees of the prime nodes $k_p$ calculated numerically (red) from the adjacency matrix using equation \eqref{3.1} exactly coincides with the degrees calculated analytically (black) from the floor function using equation \eqref{3.2} for network size $N$ = $10^4$. We observe from the figure that all prime nodes between $\fl*{\frac{N}{2}}$ and $N$ have a degree of one, all primes between $\fl*{\frac{N}{3}}$ and $\fl*{\frac{N}{2}}$ have degree as two and so on.}\label{fig:2}
\end{figure}

In figure \ref{fig:2}, we show the results for the trend in the degrees of primes for a given network size $N$ using both the equations \eqref{3.1} and \eqref{3.2}. The figure shows that the results obtained from these different methods are exactly the same. Moreover, we can see that the jumps in the values of degree at $\fl*{\frac{N}{2}}, \fl*{\frac{N}{3}}, \fl*{\frac{N}{4}}, \dots$  can be explained as coming from the nature of the floor function. We observe from figure \ref{fig:2} that all prime nodes between $\fl*{\frac{N}{2}}$ and $N$ have a degree of one, all primes between $\fl*{\frac{N}{3}}$ and $\fl*{\frac{N}{2}}$ have degree as two and so on. We also note that the jumps in the values of degrees $k_p$ can be explained in terms of the jumps in the floor function $\fl*{\frac{N}{p}}$. For example, if $N=100$, any prime between $50$ and $100$, e.g. $71$ will have $\fl*{\frac{100}{71}} = 1$, any prime between $33$ and $50$, e.g. $37$ will have $\fl*{\frac{100}{37}} = 2$, and so on.

For composite numbers the degree $k_n$ from the adjacency matrix is,
\begin{align}\label{3.3}
k_n = \sum_{m=1}^N A_{nm}
\end{align}
by definition.

In the context of composite numbers (n $>$ 1), in addition to multiples, the degree will also include the number of divisors and so we invoke the divisor function. Any composite number can be written in terms of primes as $n = \prod_{i = 1}^k p^{j_i}_i$ then the divisor function
\begin{align}\label{3.4}
    s(n) = \prod_{i = 1}^k (j_i + 1)
\end{align}
gives the number of divisors of n, which includes $n$ itself \citep{Hardy1975}. Since we do not consider self-loops, the number of divisors forming link with $n$ is $s(n) - 1$. Thus, $n$ has $\fl*{\frac{N}{n}} -1$ multiples which are less than $N$ and $s(n) - 1$ divisors, excluding $n$ itself. The degree of n, $k_n$ is the sum of the number of multiples of $n$ less than $N$ and number of divisors of $n$ excluding itself.

Thus for $n \ne 1$, we can write the degree of a node $n$ as
\begin{align}\label{3.5}
    k_n =  \fl*{\frac{N}{n}} + s(n) - 2 
\end{align}

The degree of the node $n$ = 1 is a unique case since $1$ is neither prime nor composite. The degree $k_1$ in terms of the network size is $N-1$ since 1 is connected to all the nodes in the network except itself.

\begin{figure}[htbp]
\centering
\begin{subfigure}{1\textwidth}
\centering
\includegraphics[width=1\textwidth]{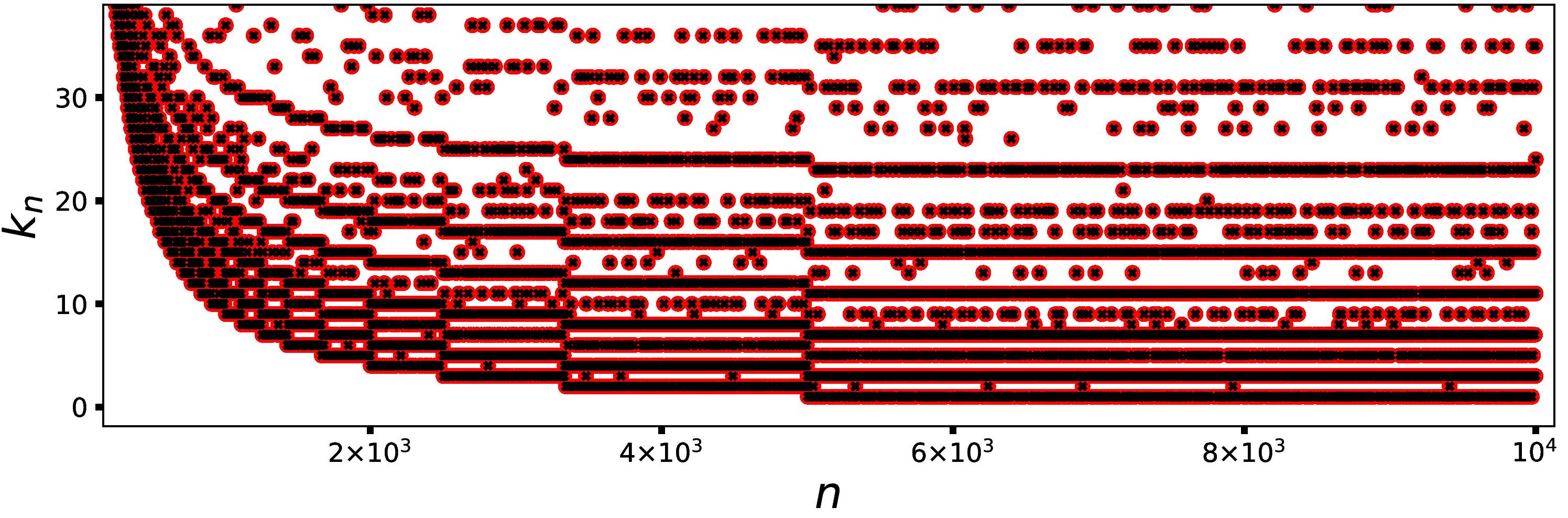}
\caption{  }\label{fig:3a}
\end{subfigure}
\centering
\begin{subfigure}{1\textwidth}
\centering
\includegraphics[width=1\textwidth]{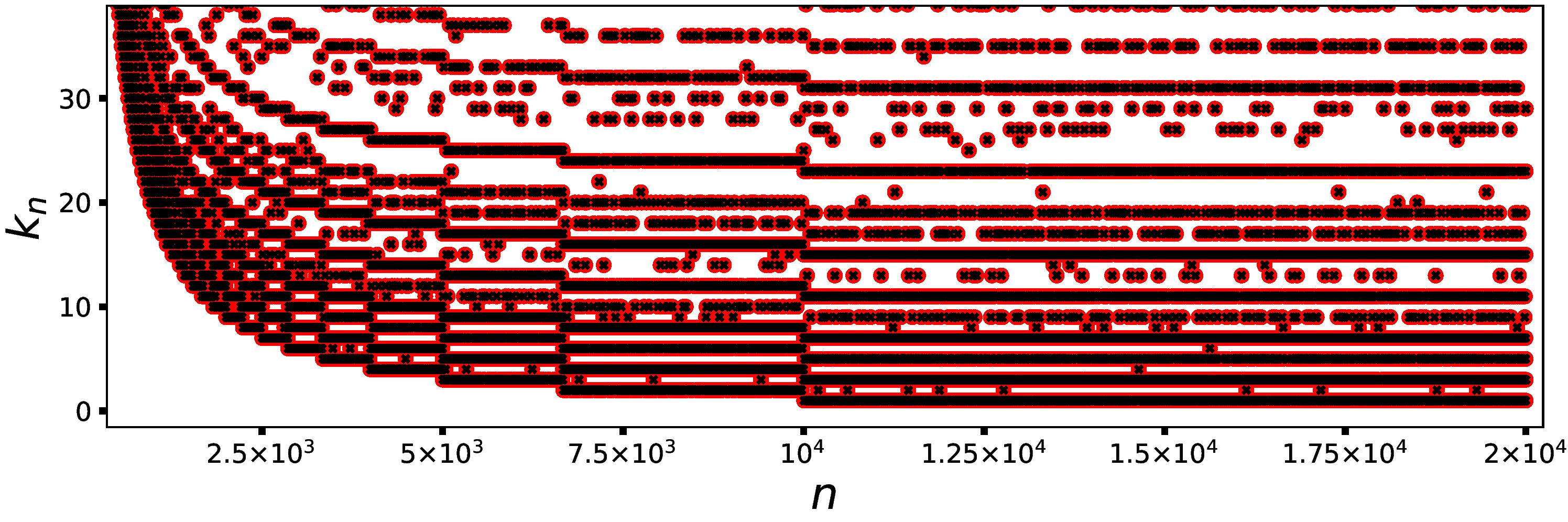}
\caption{ }\label{fig:3b}
\end{subfigure}
\caption{Degrees of nodes using the analytical (black) expression \eqref{3.5} matches exactly with the degrees obtained from adjacency matrix (red) using equation \eqref{3.3} for (a) $N$ = $10^4$, (b) $N$ = $2 \times 10^4$. We observe the whole structure stretches as N increases, giving similar patterns in both cases.}\label{fig:3}
\end{figure}

In figure \ref{fig:7}, we present the degrees of all nodes for different networks of size $N$ = $10^4$ and $N$ = $2 \times 10^4$. In addition to the jumps, we observe parallel lines and the pattern is similar for different network sizes $N$ as is clear from figures \ref{fig:3a} and \ref{fig:3b}. We note that the
jumps occur exactly at $\fl*{\frac{N}{2}}$, $\fl*{\frac{N}{3}}$, $\dots$. Comparing this figure with figure \ref{fig:2}, it is clear that the jumps come from the floor function, while the stretching formed by the bands of numbers with same degree, reflects the nature of the divisor function. This will become more evident if we plot $s(n)$ vs $n$ and $\fl*{\frac{N}{n}}$ vs $n$ and compare them with plot of $k_n$ vs $n$. We present these for a network of size $N = 10^4$ in the Appendix, which validates our conclusion.

\section{Variation of link density with the size of the network $G_N$ }
\begin{figure}[hbtp]
\centering
\begin{subfigure}{1\textwidth}
\includegraphics[width=1\textwidth]{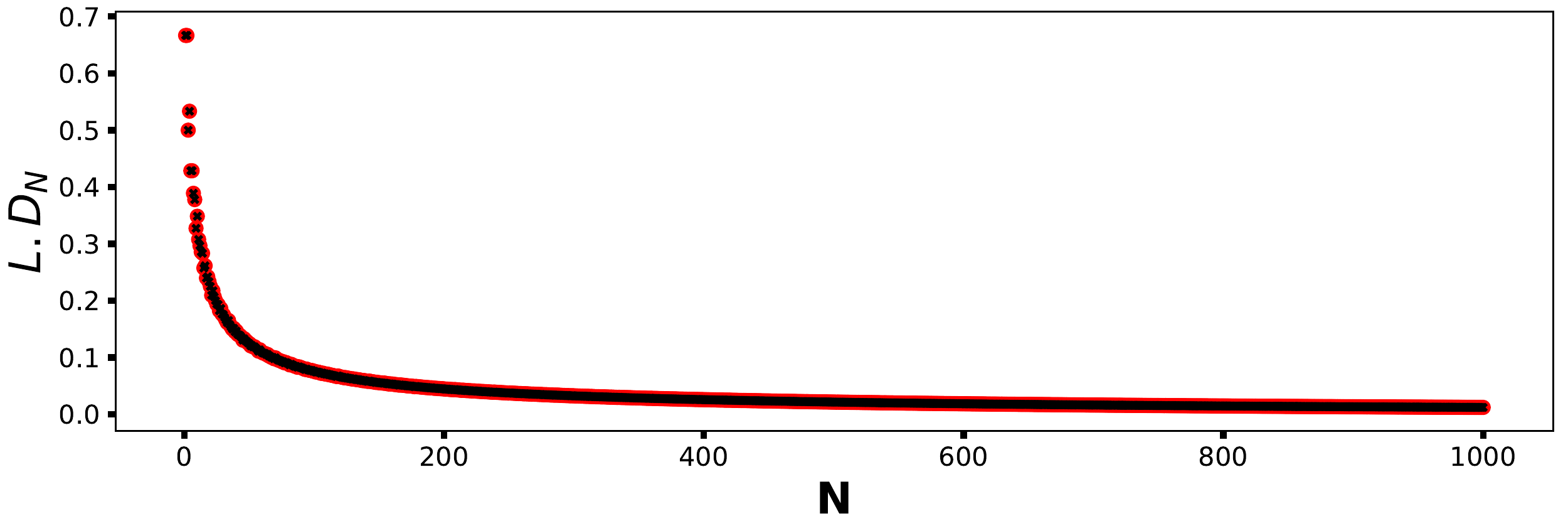}
\caption{  }\label{fig:4a}
\end{subfigure}
\vspace*{\fill} 
\begin{subfigure}{1\textwidth}
\includegraphics[width=1\textwidth]{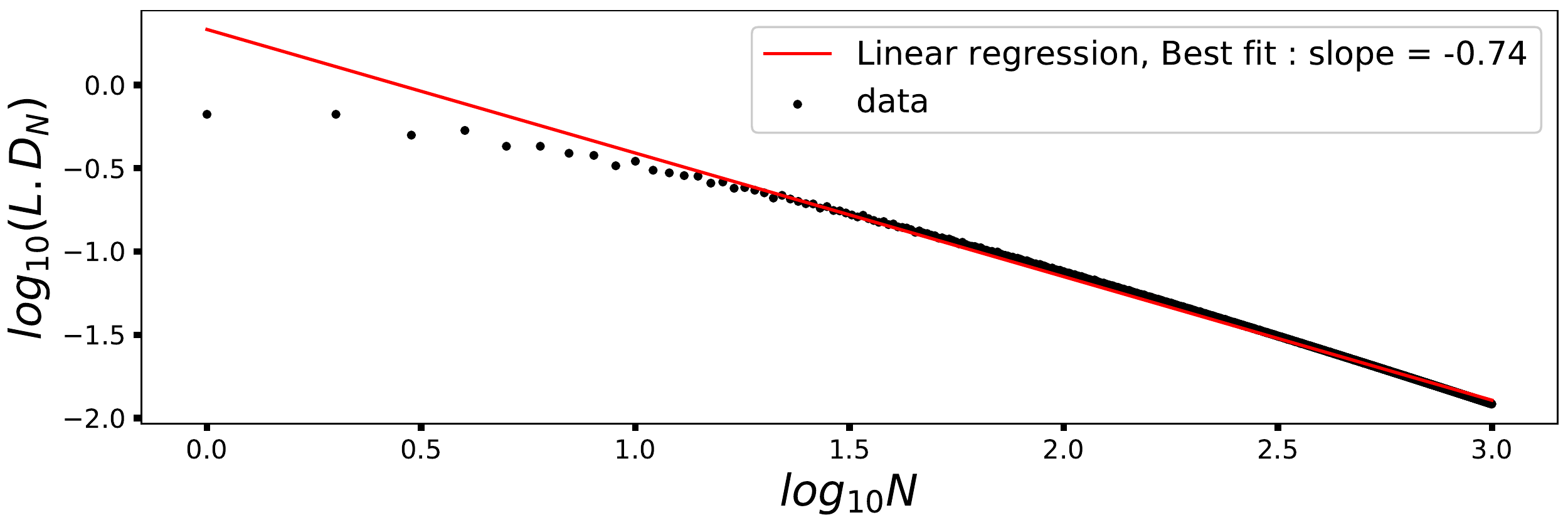}
\caption{ }\label{fig:4b}
\end{subfigure}
\caption{ (a) There is an exact match of change in link density with the size of the network obtained from the adjacency matrix (red) using equation \eqref{4.1} with the values obtained using the analytical expression (black) in equation \eqref{4.2}. (b) The log-log plot of variation in link Density with the size of the network. The slope of the straight-line portion is approximately -3/4, which indicates the scaling of link density with size of network.}\label{fig:4}
\end{figure}

One of the fundamental global characteristics of a network is its link density. It is defined as the ratio of the number of links present in a network to the number of possible links between all node pairs n,m in the network of size $N$ \citep{Newman2010}. Thus from adjacency matrix, we can write 

\begin{align}\label{4.1}
\text{Link density, } L.D_N = \sum_{n,m}\frac{0.5A_{nm}}{^NC_2}
\end{align}

Following the discussions in the above section, for $G_N$, the total number of links in $G_N$ is the sum of the number of multiples less than $N$ (which is $\fl*{\frac{N}{i}}-1$) of all the nodes in $G_N$. Thus we get,     
\begin{align}\label{4.2}
\text{Link density, } L.D_N   =  \sum_{i=1}^N \frac{ \fl*{\frac{N}{i}} } {(^NC_2)}  - \frac{2}{(N-1)} .
\end{align}

In figure \ref{fig:4a}, we plot the link density as a function of network size, using both the above expressions. It is clear that the agreement is good. From the log-log plot shown in figure \ref{fig:4b}, we get the slope of the straight-line portion as -0.74 or approximately -3/4. Thus we can say for large $N$, the link density scales with $N$ as,

\begin{align}\label{4.3}
LD(N) \approx  N^{-3/4}
\end{align}

This would mean as $N$ increases, link density goes to zero following a power law. We note this is different from the scaling index of -1 reported in real-world networks \citep{Blagus2012}.

\section{Local clustering coefficient and stretching similarity}

The local clustering coefficient of any node, $n$ is the ratio of the number of neighboring pairs of $n$ which are connected among themselves ($e_n)$ to the number of possible links among the neighbors of $n$ \citep{Newman2010}. In terms of elements of the adjacency matrix, the local clustering coefficient of any node $n$ is given by

\begin{align}\label{5.1}
c_n = \frac{e_n}{^{k_n}C_2} = \sum_{s \ne t = 1}^N \frac{A_{ns}A_{st}A_{tn}}{{}^{k_n}C_2} 
\end{align}

We present below the derivation of analytical expression for the local clustering coefficient of any node $n$ is discussed below:

We start by considering the number of links among the divisors of $n$ (call it $N_{s(n)}$). There are $s(n)-1$ divisors of $n$ excluding $n$ itself. Let $m$, $m'$ be two divisors of $n$ ($= \prod_{i = 1}^k p^{j_i}_i$). Then, we can write $m = \prod_{\alpha = 1}^k p^{l_\alpha}_\alpha$ and $m' = \prod_{\alpha = 1}^k p^{l_\alpha'}_\alpha$ where $0 \le l_\alpha, l_\alpha' \le j_\alpha $. If $m'$ is a divisor of $m$ (denoted $m'|m$) and there is an edge between $m$ and $m'$, then $l_\alpha \le l_\alpha'$ for each $\alpha = 1, \dots, k$. For fixed $m$, the number of divisors of $m$ is $s(m)-1$ = $\prod_{\alpha=1}^k (l_\alpha + 1) - 1$ . Summing this over all the divisors $m$ of $n$ gives, $\displaystyle{\sum_{m|n} (s(m) - 1)}$ $=\displaystyle{\sum_{m|n} s(m) - s(n)}$. But this includes the case when $m=m'$ for each divisor $m$. So, subtracting the number of such $m'$ from this gives the number of edges between the 
divisors of n, 
\begin{align}\label{5.2}
N_{s(n)} = \sum_{m|n} s(m) - s(n) - (s(n)-1)
\end{align}
\begin{align}\label{5.3}
  = \left( \sum_{l_1 = 0}^{j_1}\sum_{l_2 = 0}^{j_2} \dots \sum_{l_k = 0}^{j_k} \left[\prod_{\alpha=1}^k (l_\alpha + 1) \right] \right) - (2s(n)-1).
\end{align}
\begin{align}\label{5.31}
  = \frac{\displaystyle{\prod_{i=1}^k} \  (j_i+1)(j_i+2)}{2^k} - 2s(n)+1.
\end{align}

Now, we look at the number of links among the multiples of $n$ (call it $N_{M(n)}$) Here $M(n) = \fl*{\frac{N}{n}}$ is the number of multiples of $n$ which are less than $N$, including $n$ itself. The multiples of $n$ less than $N$ excluding $n$ are $2.n, \  3.n, \  \dots, \  M(n).n$.  There is an edge between $\alpha.n$ and $\beta.n$ if either $\alpha$ divides $\beta$ or $\beta$ divides $\alpha$. To find the number of such edges (called $N_{M(n)}$), let us look at a divisibility graph $G_N$ with nodes $\{1,\dots,M(n)\}$  (call this graph $G_{M(n)}$). Observe that the number of such edges $N_{M(n)}$ is the number of edges in the graph $G_{M(n)}$ after removing all the edges connected to $1$. 

Since each vertex $j$ in $G_{M(n)}$ is connected to $\fl*{\frac{M(n)}{j}}-1$ number of multiples, $\displaystyle{\sum_{j=1}^{M(n)} \left(\fl*{\frac{M(n)}{j}}-1\right)}$ gives the total number of edges in $G_{M(n)}$. Removing the vertex $1$, 
\begin{align}\label{5.4}N_{M(n)} = \sum_{j=2}^{M(n)} \left(\fl*{\frac{M(n)}{j}}-1\right) = \sum_{j=2}^{M(n)} \fl*{\frac{M(n)}{j}} - (M(n) - 1) \end{align}

Next, we discuss the number of links between the divisors of $n$ and the multiples of $n$ ( call it, $N_{int}$ ). Each divisor of $n$ is also a divisor of any multiple of $n$. So, there is an edge connecting each divisor of $n$ to each of its multiple. Hence, the number of edges connecting a multiple of $n$ to a divisor of $n$ is given by
\begin{align}\label{5.5}
N_{int} = (M(n)-1)(s(n)-1)
\end{align}

Finally, the number of links among the neighbors of $n$ ( $e_n$ ), is the sum of the number of links among the divisors of $n$ ( $N_{s(n)}$ ), the number of links among the multiples of $n$ ( $N_{M(n)}$ ) and the number of links connecting a multiple of $n$ to a divisor of $n$ ( $N_{int}$ ).
\begin{align}\label{5.6}
e_n = N_{s(n)} + N_{M(n)} + N_{int}
\end{align}
\begin{align}\label{5.7}
e_n = \sum_{m|n} s(m) - 2s(n)+1 + \sum_{j=2}^{M(n)} \fl*{\frac{M(n)}{j}} + (M(n)-1)(s(n)-2)
\end{align}
for any node $n \ ( = \prod_{i = 1}^k p^{j_i}_i )$, in the $G_N$ in terms of the divisor function and floor function. Here 
\begin{align}\label{5.8}
\sum_{m|n} \  s(m) = \frac{\displaystyle{\prod_{i=1}^k} \  (j_i+1)(j_i+2)}{2^k}.
\end{align}

\begin{figure}[htbp]
\centering
\begin{subfigure}{1\textwidth}
\centering
\includegraphics[width=1\textwidth]{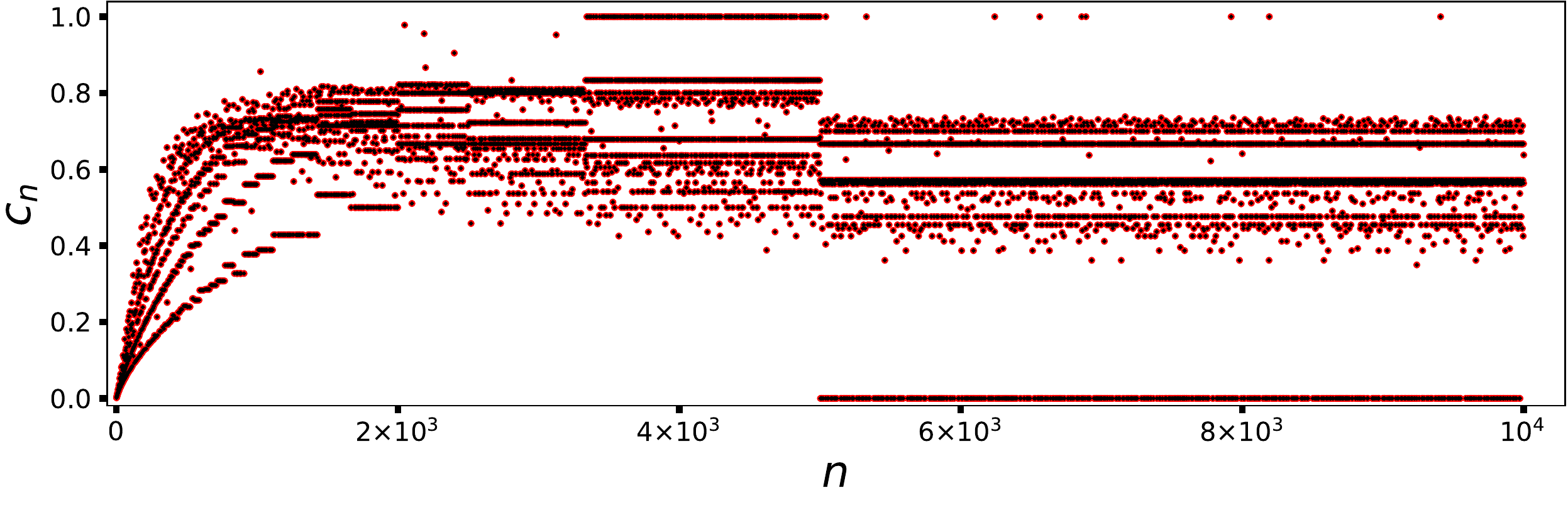}
\caption{  }\label{fig:5a}
\end{subfigure}
\centering
\begin{subfigure}{1\textwidth}
\centering
\includegraphics[width=1\textwidth]{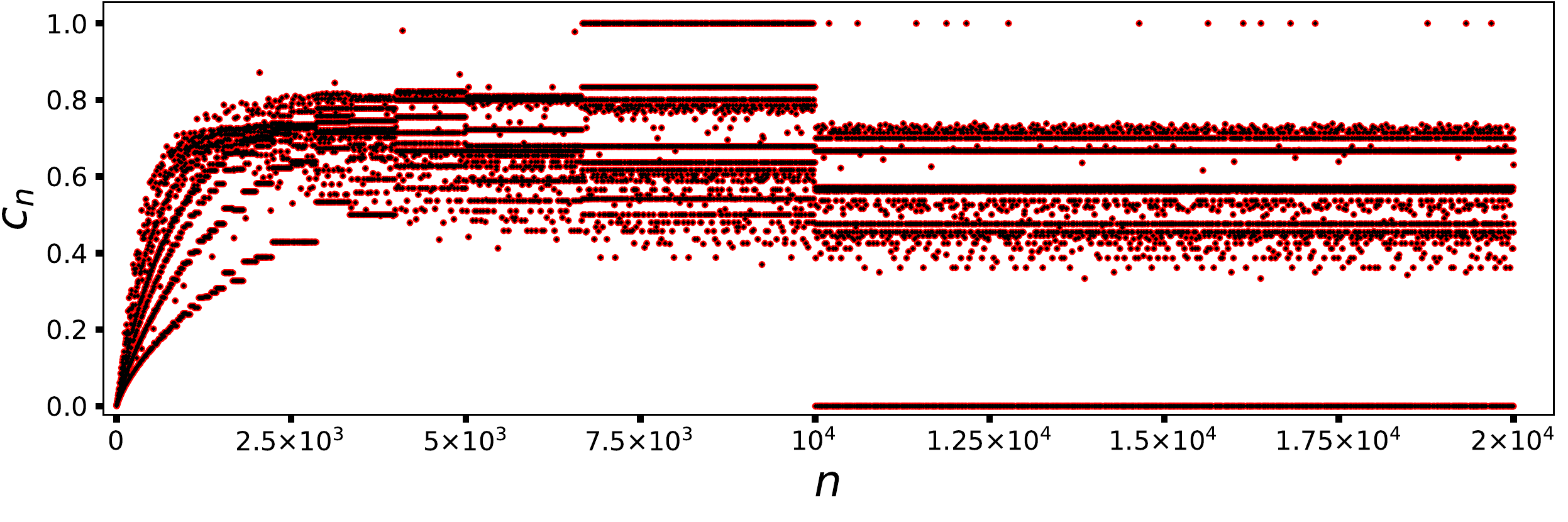}
\caption{ }\label{fig:5b}
\end{subfigure}
\caption{Local clustering coefficients $c_n$ of nodes $n$ in the divisibility network of size (a) $N$ = $10^4$ and  (b) $N$ = $2 \times 10^4$. It is clear from (a), (b) that the values obtained from the analytical (black) equation \eqref{5.9} coincides exactly with the $c_n$ calculated from adjacency matrix (red) using equation \eqref{5.1}. We note the stretching similarity in the distributions of $c_n$ as is clear from (a) and (b).}\label{fig:5}
\end{figure}

Hence, the local clustering coefficient of any node $n$ is analytically derived to be
\begin{align}\label{5.9}
c_n = \frac{N_{s(n)} + \displaystyle{\sum_{j=2}^{M(n)} \fl*{\frac{M(n)}{j}}} + (M(n)-1)(s(n)-2)}{^{k_n}C_2}
\end{align}
where 
\begin{align}\label{5.10}
N_{s(n)} = \sum_{m|n} s(m) - 2s(n)+1 \text{ and } M(n) = \fl*{\frac{N}{n}}.
\end{align}
The above formula is computationally verified in figure \ref{fig:5} where we show the coincidence of $c_n$ computed using the expression derived by us and also using the standard expression from the adjacency matrix (equation \eqref{5.1}). 

A new kind of similarity called stretching similarity was reported in the local clustering coefficient of $G_N$ \citep{Shekatkar2015}. We find this can be explained using our analytic expression. The expression for $c_n$ contains the divisor function $s(n)$, the number of multiples of $n$ including $n$, ($M(n)$) and the number of divisors of the divisors of $n$.

The main features observed in the stretching similarity in figure \ref{fig:5},are jumps, stretching and an overall irregularity. We observe sharp jumps at $\fl*{\frac{N}{2}}$, $\fl*{\frac{N}{3}}$, $\fl*{\frac{N}{4}}$, $\dots$ that comes from the floor function. In addition there are bands that comes from the number of divisors $s(n)$, which depends only on the powers of the primes in the prime factorisation of a given node $n$. The closer nodes with the same number of divisors and prime powers lie in the same band. As N increases the pattern stretches giving rise to similarity which is seen in the degrees of nodes also. As mentioned earlier the plot of divisor function vs nodes also shows similar bands (see Appendix). 

The overall irregularity is caused by nodes which have the same number of divisors but different prime powers, for example when $N=50$, the nodes $32$ and $45$ have the same number of divisors $s(32) = 6, s(45) = 6$. But they have different prime powers, since $32 = 2^5$ and $45=3^25$. Hence, their local clustering coefficients are different, $c_{32} = \frac{20}{21} = 0.9523 $ and $c_{45} = \frac{47}{66} = 0.7121$.
\begin{figure}[htbp]
\centering
\begin{subfigure}{1\textwidth}
\centering
\includegraphics[width=1\textwidth]{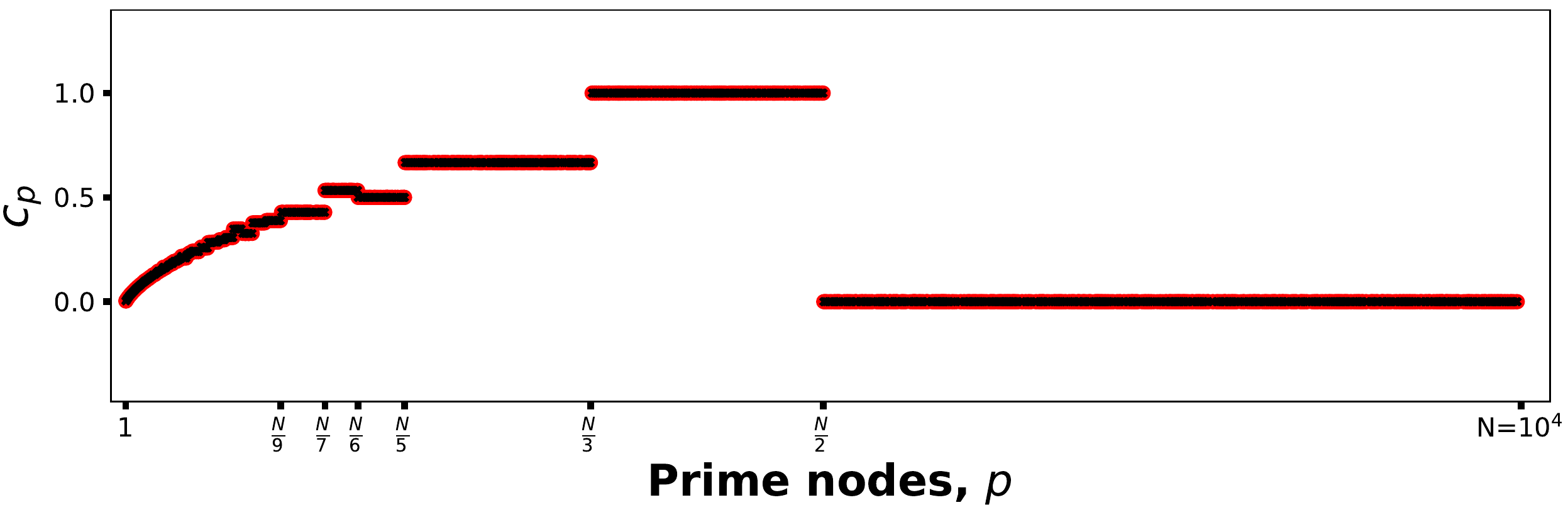}
\caption{  }\label{fig:6a}
\end{subfigure}
\centering
\begin{subfigure}{1\textwidth}
\centering
\includegraphics[width=1\textwidth]{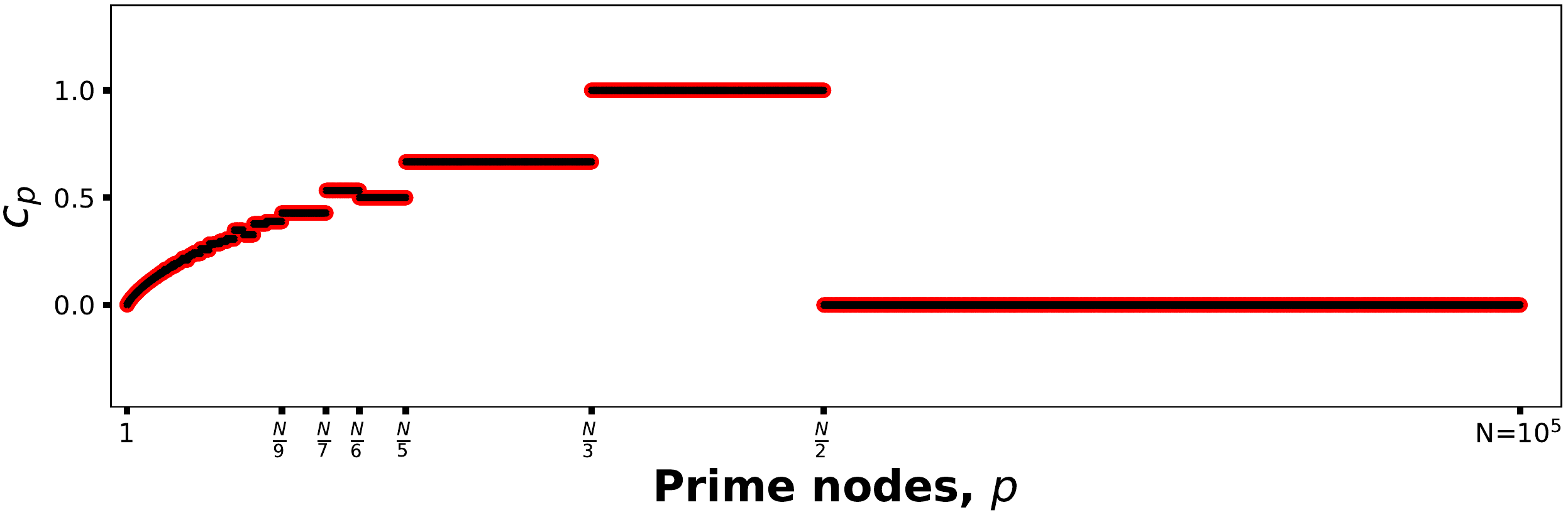}
\caption{ }\label{fig:6b}
\end{subfigure}
\caption{(a) Local clustering coefficients of the prime nodes $c_p$ calculated numerically (red) from the adjacency matrix using equation \eqref{5.1} exactly coincides with the degrees calculated analytically (black) from the floor function using equation \eqref{05.12} for a network size of (a) $N$ = $10^4$ and (b) $N$ = $10^5$. We observe a trend from the figure that all primes $p$ between $\fl*{\frac{N}{2}}$ and $N$, have $c_p=0$, all primes between $\fl*{\frac{N}{3}}$ and $\fl*{\frac{N}{2}}$, have $c_p=1$ and all primes $p$ between $\fl*{\frac{N}{5}}$ and  $\fl*{\frac{N}{3}}$, have $c_p = \frac{2}{3}$, and so on.}\label{fig:6}
\end{figure}

Next we discuss the particular case of clustering coefficients of prime numbers.

For a given $N$, if p is a prime such that $\fl*{\frac{N}{2}} < p \le N$, then $c_p = 0$. This is because 1 is its only neighbor, making $e_p=0$. Similarly, if p is a prime such that $\fl*{\frac{N}{3}} < p \le \fl*{\frac{N}{2}}$, then $c_p = 1$. This is because such a prime's neighbors in the network are 1 and 2p which are connected in the network making $c_p=1$.

For p such that $\fl*{\frac{N}{a+1}} < p \le \fl*{\frac{N}{a}}$, $M(p) = a$. Then, using $s(p) = 2$ and $M(p) = a$ in equation \eqref{5.7}, we get $e_p$ = $\displaystyle{\sum^{a}_{j=2} \fl*{\frac{a}{j}}}$. This gives $c_p$ as 
\begin{align}\label{05.12}
    c_p = \frac{\displaystyle{\sum^{a}_{j=2} } \fl*{\frac{a}{j}}}{{}^aC_2}
\end{align}
since $k_p = \fl*{\frac{N}{p}} = a$. Figure \ref{fig:6} computationally verifies the above equation for different network sizes of $N$ = $10^4$ and $N$ = $2 \times 10^4$. This figure also shows other trends in the prime numbers lying between different intervals. For example, for $N=1000$, all primes between $200$ and $333$ have $c_p = \frac{2}{3}$, all primes between $166$ and $200$ have $c_p = \frac{1}{2}$, and so on.

\begin{figure}[htbp]
\centering
\includegraphics[width=1\textwidth]{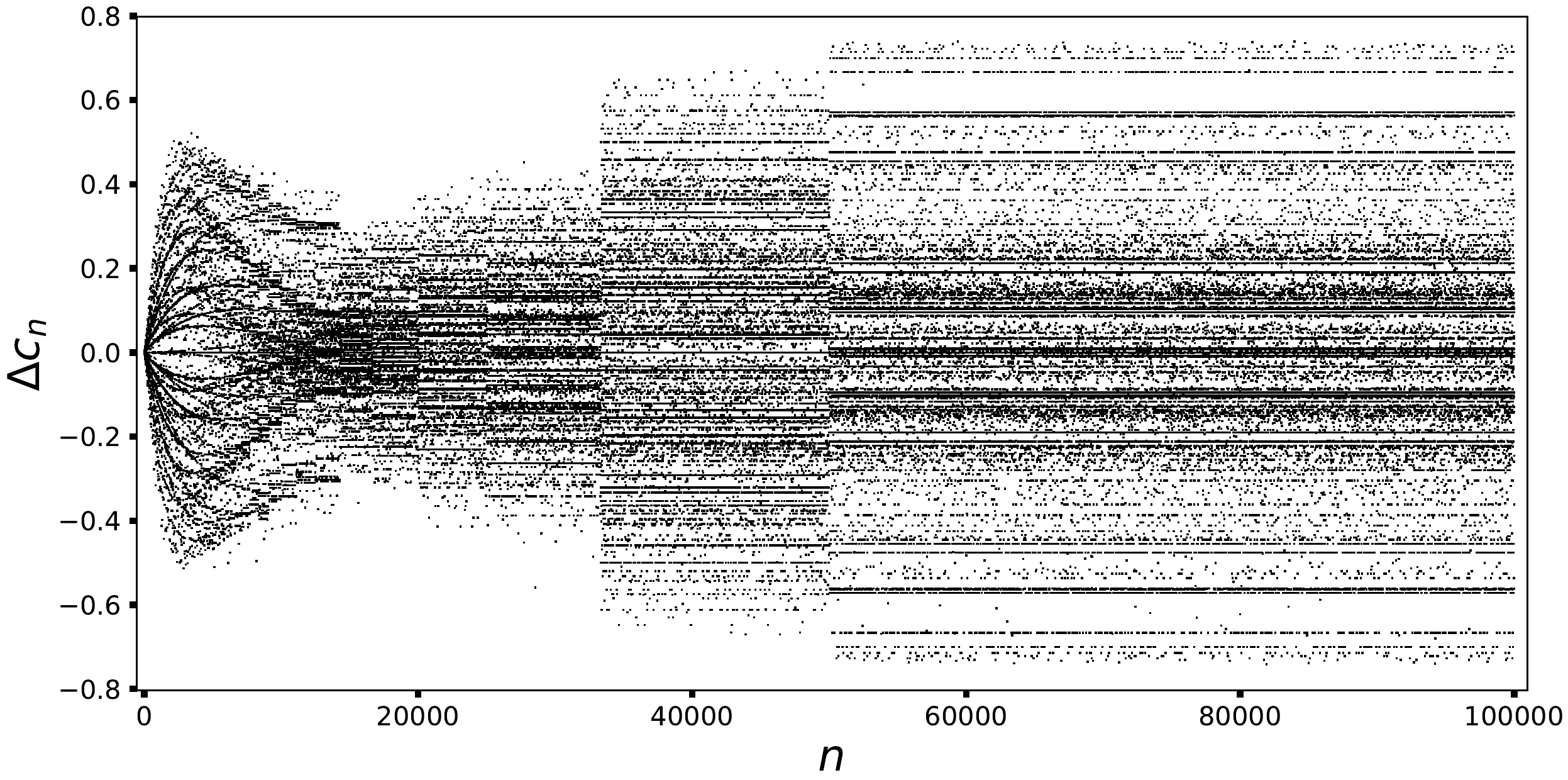}
\caption{Plot of $\Delta c_n =  c_n - c_{n+1}$ vs $n$ for a network size of $10^5$. The statistical symmetry reported in \citep{Shekatkar2015} is seen in this figure, plotted using the analytic expressions.}\label{fig:7}
\end{figure}

Now, we discuss a type of symmetry reported in the difference between the local clustering coefficients of consecutive numbers, called the stretching symmetry, in terms of our derived expression \eqref{5.9}.
Figure \ref{fig:7}, gives the plot of $\Delta c_n =  c_n - c_{n+1}$ vs $n$. Ref \cite{Shekatkar2015} numerically showed that there is stretching similarity in this plot. They also observed that this figure is statistically symmetric about the line $\Delta c_n = 0$. We give an explanation for the reason behind this statistical symmetry and stretching similarity in $\Delta c_n$, based on our analytical expression. 

In the $c_n$ vs $n$ plot, there is only a finite number of bands between a given $\fl*{\frac{N}{a+1}}$ and $\fl*{\frac{N}{(a)}}$ ($a=1,2,3,\dots$). If $c_n$ and $c_{n+1}$ lie in the same band, then $\Delta c_n =  0$. For example, let us consider a network size of $N=100$ (Figure \ref{fig:6}), then $93$ \& $94$ lie between $\fl*{\frac{N}{2}}$ and $\fl*{\frac{N}{1}}$ and have the same $c_n = \frac{2}{3}$. Therefore $\Delta c_{93} = 0$. Other such pairs are $94$ \& $95$, $86$ \& $87$ and $85$ \& $86$ which lie in the same band $c_n = \frac{2}{3}$, resulting in $\Delta c_n = 0$.  Similarly, there are many such pairs depending on the size of the network.

On the other hand, if $c_n$ lies on a band above(below) $c_{n+1}$ , then $\Delta c_n$ is positive (negative) (examples are, $c_{82} = \frac{2}{3}$, $c_{83} = 0$, so $\Delta c_{82} = \frac{2}{3}$ and $c_{73} = 0$  $c_{74} = \frac{2}{3}$, so $\Delta c_{73} = -\frac{2}{3}$). Moreover, as seen from figure \ref{fig:5}, for a fixed $a$, $|\Delta c_n|$ takes only a finite number of values and we also found that as $N$ changes the count of numbers with $\Delta c_n = 0$ changes in a systematic way. This finiteness is again due to the nature of the divisor function, which can take only finite number of values for $n$ less than $N$.

From our analytical expressions we can obtain $\Delta c_n = 0$ as, \begin{align}\begin{split}\label{5.11}\Delta c_n = \frac{\displaystyle{\sum_{m|n} s(m) } - 2s(n) + \sum_{j=2}^{M(n)} \fl*{\frac{M(n)}{j}} + (M(n)-1)(s(n)-2)}{0.5k_n(k_n-1)} \\ - \   \frac{\displaystyle{\sum_{m'|n+1} s(m')} - 2s(n+1) + \sum_{j=2}^{M(n+1)} \fl*{\frac{M(n+1)}{j}} + (M(n+1)-1)(s(n+1)-2)}{0.5k_{n+1}(k_{n+1}-1)} \end{split}\end{align} Let us take $n,n+1$ such that $\fl*{\frac{N}{a+1}} < n < n+1 \le \fl*{\frac{N}{a}}$ for some positive integer $a$, then $M(n) = M(n+1) = a$. Then, simplifying the above expression, we get \begin{align}\begin{split}\label{5.12}\Delta c_n = \frac{\left(\displaystyle{\prod_{i=1}^k} \  (j_i+1)(j_i+2)/2^k \right) - 2s(n) + (a-1)(s(n)-2)}{0.5(a + s(n) - 2)(a + s(n) - 3)}\\ - \  \frac{\left(\displaystyle{\prod_{i=1}^{k'} \  (j_i'+1)(j_i'+2)/2^{k'}}\right) - 2s(n+1) + (a-1)(s(n+1)-2)}{0.5(a + s(n+1) - 2)(a + s(n+1) - 3)} \end{split} \end{align} From the above equation, we see that $\Delta c_n = 0$ when \begin{align}\label{5.13}
\frac{\displaystyle{\prod_{i=1}^k} \  (j_i+1)(j_i+2)}{2^k}  = \frac{\displaystyle{\prod_{i=1}^{k'}} \  (j_i'+1)(j_i'+2)}{2^{k'}} \end{align} and \begin{align}\label{5.14} s(n) = s(n+1)\end{align} 
i.e, $n,n+1$ have the same prime powers and the same number of divisors. 

This supplements the result shown numerically in Ref. \citep{Shekatkar2015} that, there are many consecutive numbers $n,n+1$ such that $\Delta c_n = 0$ which means that these numbers have the same prime powers and the same number of divisors and lie in the same interval. 

We relate this property to the theorem proved by Heath-Brown (1984), that $s(n) = s(n+1)$ infinitely many times \citep{Heath1984}. Some examples are $s(93)=s(94)=4$, $s(94)=s(95)=4$, $s(86)=s(87)=4$, $\dots$. This means many consecutive natural numbers can have the same $c_n$, if their prime powers are also equal, giving $\Delta c_n = 0$. Examples are $\Delta c_{93} = 0$, $\Delta c_{94} = 0$, $\Delta c_{86} = 0$,$\dots$. Further we observe that there are different bands with $ \Delta c_n \ne 0$, in the $\Delta c_n =  c_n - c_{n+1}$ vs $n$ plot. This could mean that there are also infinitely many numbers with $s(n) = s(n+1)+k$ for some $k$ other than zero and this can be considered as an extension of the Heath-Brown theorem. For example, for a network of size $N=100$, $\Delta c_{82} = \frac{2}{3}$ and $s(82)=4$, $s(83)=2$, so in this case $k=2$, since $s(82) = s(83)+2$. Another example, is for $N=5$, $\Delta c_{3} = -1$ and $s(3)=2$, $s(4)=3$, so in this case $k=-1$, since $s(3) = s(4)-1$. 

We find $ \Delta s(n)$ also shows a symmetry about $\Delta s(n) = 0$ as shown in Appendix.

\section{Betweenness Centrality from adjacency matrix for $G_N$:}

The approach presented here can, in principle, be further extended to derive expressions for any of the characteristic measures of $G_N$. In this section, we consider one of them, namely, the betweenness centrality. In general, it is computationally difficult to calculate the number of geodesic paths and find the betweenness centrality $x_n$ of each node $n$ in any network and most often approximation algorithms are used to compute the same. However, for the divisibility number network $G_N$, it being a deterministic network, we can easily compute $x_n$ from the adjacency matrix itself. By definition,
\begin{align}\label{6.1}
x_n = \sum_{s,t \ne n} \frac{n_{st}^n}{g_{st}} \ \ \ \frac{1}{(N-1)(N-2)}
\end{align}
where $g_{st}$ is the number of geodesic paths between s and t and $n_{st}^n$ of them pass through $n$ \cite{Newman2010}. These two can be calculated directly from the adjacency matrix, exploiting the fact that we know the geodesic path between any two nodes in $G_N$.

As mentioned, $g_{st}$ is the number of geodesic paths between $s$ and $t$ and $n_{st}^n$ is the number of geodesic paths between $s$ and $t$ which pass through $n \ne s,t$. If $s$ and $t$ are connected, then the geodesic path is of length one and it cannot pass through another $n$. So, $n_{st}^n$ = 0, if $s$ and $t$ are connected. If $s$ and $t$ are not connected then, geodesic path is of length 2, since there is a path between any two nodes in $G_N$ via the node $1$. The number of paths of length 2 between $s$ and $t$ is given by $[A^2]_{st}$  \citep{Newman2010}. Hence, $g_{st}$ = $[A^2]_{st}$, if $s$ and $t$ are not connected. And $n_{st}^n$ can be written in terms of the adjacency matrix as $n_{st}^n$ = $A_{sn}A_{nt}$, which is one only if there is a path of length two via $n$ between $s$ and $t$.  

\begin{figure}[hbtp]
\centering
\includegraphics[width=1\textwidth]{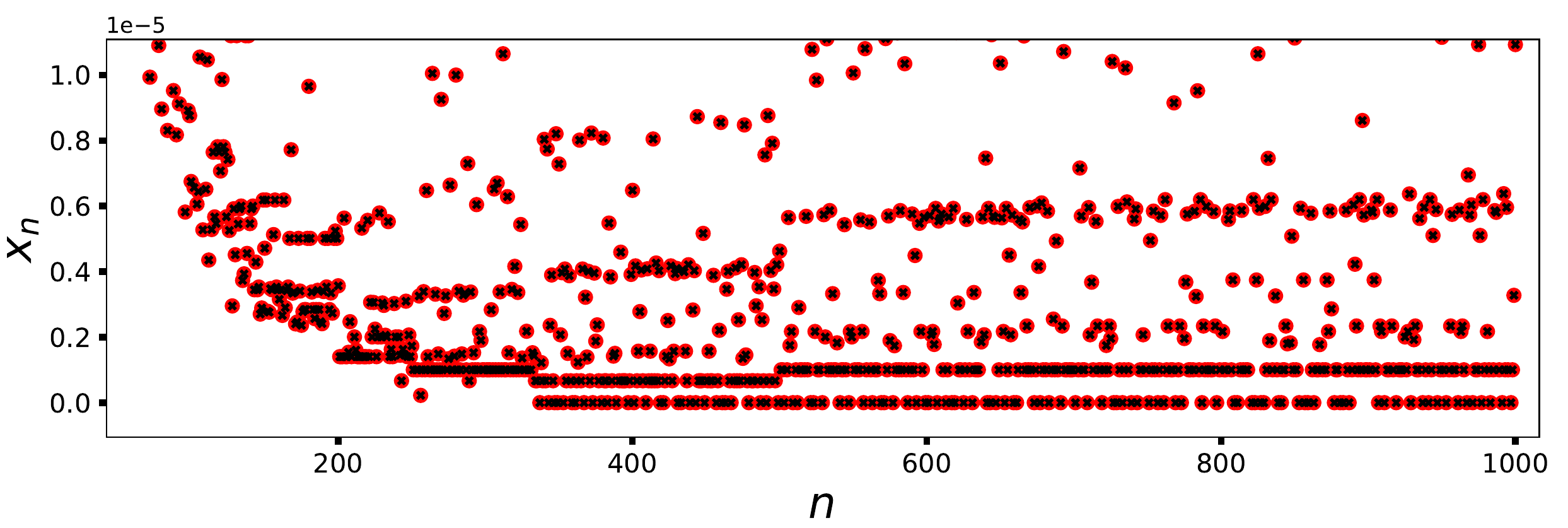}
\caption{Betweenness Centrality $x_n$ vs nodes $n$ for the network of size $N=1000$. The values calculated from the adjacency matrix (black) using equation \eqref{6.2} coincide exactly with the values calculated from python-networkx function (red) using Ulrik Brandes Algorithm \citep{Brandes2001}.}\label{fig:8}
\end{figure}

Hence, we can simply write $x_n$ as
\begin{align}\label{6.2}
    x_n = \sum_{s,t \ne n} \frac{(1-A_{st})A_{sn}A_{nt}}{[A^2]_{st}} \ \ \ \frac{1}{(N-1)(N-2)}
\end{align}
The $(1-A_{st})$ part ensures that this term is zero, when $s$ and $t$ are connected. Equation \eqref{6.2} gives a way to compute $x_n$ for the divisibility graph of natural numbers directly from the adjacency matrix (Figure \ref{fig:8}).

\section{Conclusion}
We present a study on the characteristic measures of the natural number network and illustrate how these measures can reveal the trend in the behaviour of natural numbers. The analytic expressions derived for the network measures like degree, clustering coefficient and link density are dependent the floor functions and divisor functions of natural numbers. Thus our study deepens the connection between number theory and graph theory, and we could relate the measures of $G_N$ to properties of natural numbers. We present the trends in the measures specifically for prime numbers, and in general for composites. Our results are validated using standard methods using elements of the adjacency matrix of the network. 

We could show how the link density of the network scales with the size of the network for large $N$. Our analysis also could explain the stretching similarity reported in previous work for local clustering coefficients of the nodes as $N$ is increased. We identify the same type of similarity in the degrees of nodes also and relate that to the properties of the divisor functions and floor functions of natural numbers. We could also see a possible extension of the theorem by Heath-Brown from the plot of the difference in clustering coefficients of two consecutive numbers, as there might also be infinitely many numbers with s(n) = s(n+1)+k for some k other than 1. We foresee further connections between number theory and graph theory along similar lines.

\section*{Acknowledgment}
AR would like to thank Department of Science and Technology, Government of India for INSPIRE scholarship.

\appendix

\section*{Appendix}

We present below the trends in the divisor functions $s(n)$ and floor functions $\fl*{\frac{N}{n}}$ with $n$. The trends presented in the measures of the network $ G_N$, obtained from their analytic expressions, are similar. Thus, this further supports our inference that the trends in the measures of the network $G_N$, can be related to the properties of $s(n)$ and $\fl*{\frac{N}{n}}$.

\begin{figure}[htbp]
\centering
\begin{subfigure}{1\textwidth}
\includegraphics[width=1\textwidth]{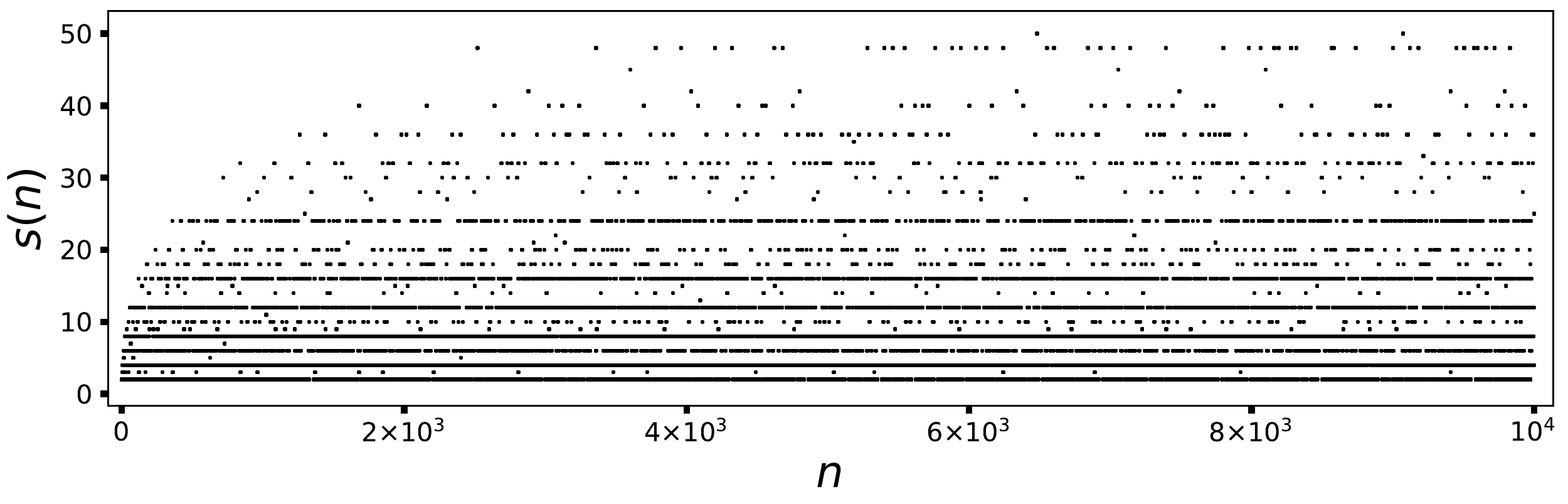} 
\caption{  }\label{fig:9a}
\end{subfigure}
\vspace*{\fill} 
\begin{subfigure}{1\textwidth}
\includegraphics[width=1\textwidth]{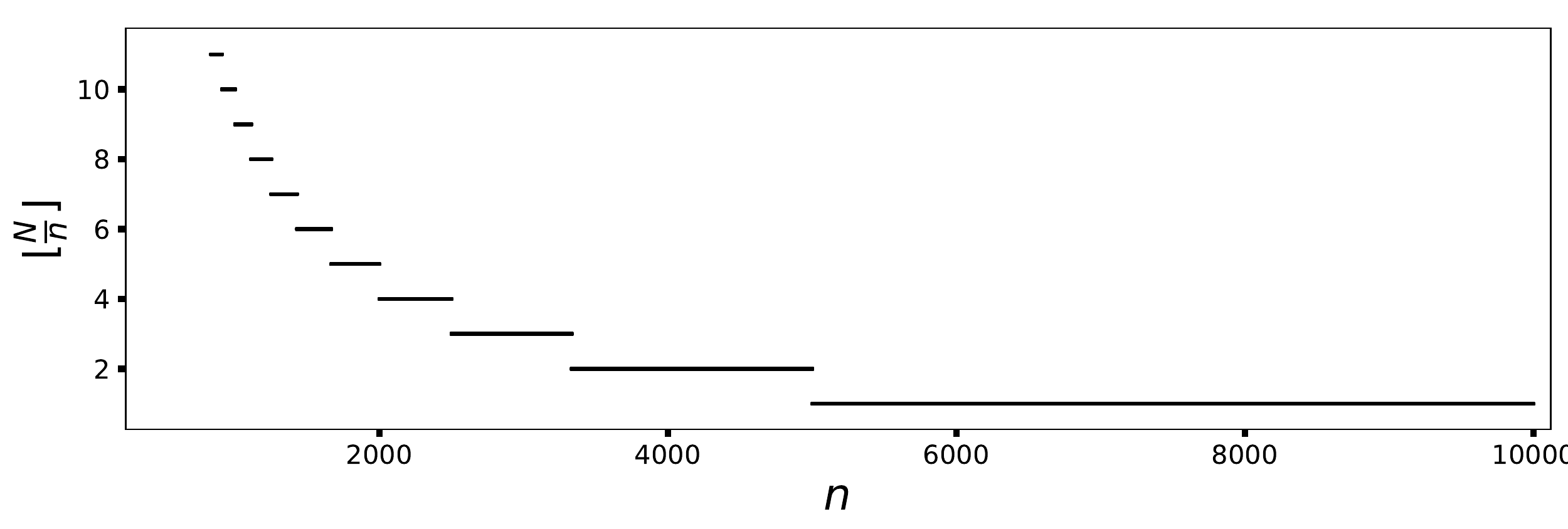}
\caption{  }\label{fig:9b}
\end{subfigure}
\vspace*{\fill} 
\begin{subfigure}{1\textwidth}
\includegraphics[width=1\textwidth]{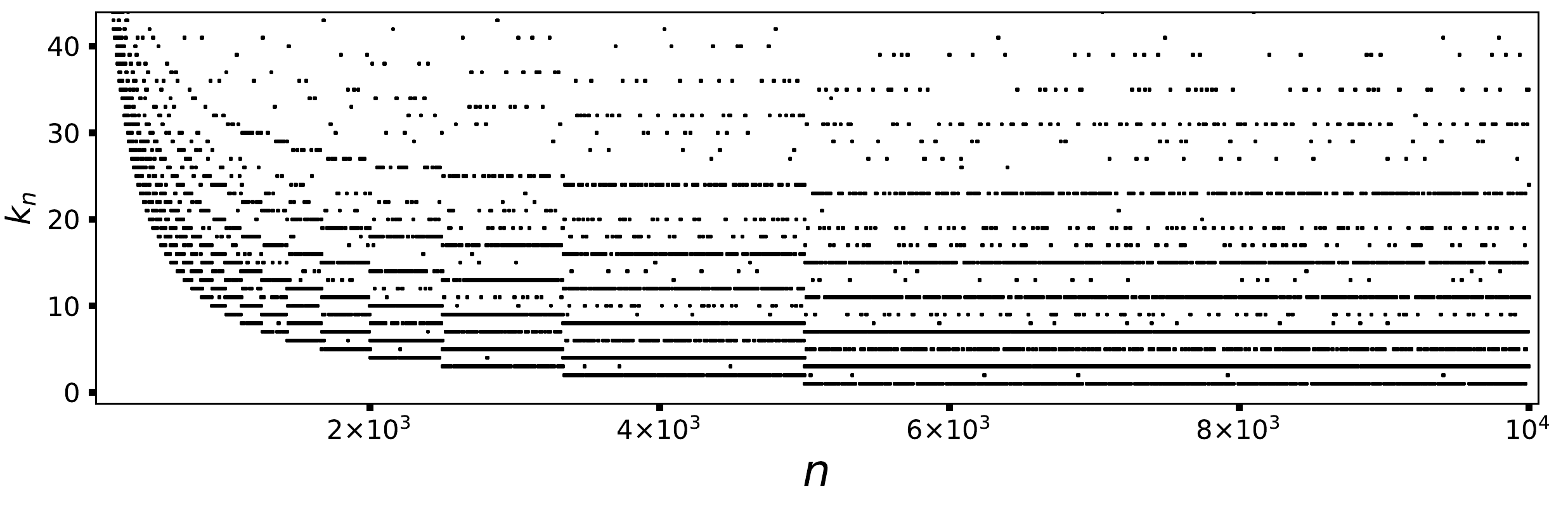}
\caption{  }\label{fig:9c}
\end{subfigure}
\caption{(a) Plot of divisor function, $s(n)$ vs nodes, $n$ for $N=10^4$. We observe that there are parallel lines in the plot. (b) Plot of $\fl*{\frac{N}{n}}$ vs nodes, $n$ for $N=10^4$, has jumps. (c) Plot of degree, $k_n = s(n) + \fl*{\frac{N}{n}} - 2$ for $N=10^4$. This has both the jumps and parallel bands derived from the above two plots.}\label{fig:9}
\begin{flushleft}
We can see from figures \ref{fig:9a}, \ref{fig:9b} that there are no characteristic jumps in the divisor function but only in the floor function. So the jumps and overall decrease in the degree arise due to the properties of floor function while the overall pattern comes from the divisor function.
\end{flushleft}
\end{figure}

\begin{figure}
\begin{flushleft}
We also present the plot of $ \Delta c_n = c_n - c_{n+1}$  vs nodes $n$ for the network of size $N=10^5$ and compare that with the plot of $\Delta s(n) = s(n) - s(n+1)$ vs $n$. It is clear from the figure that the statistical symmetry comes from the nature of the divisor function.
\end{flushleft}
\centering
\begin{subfigure}{1\textwidth}
\includegraphics[width=1\textwidth]{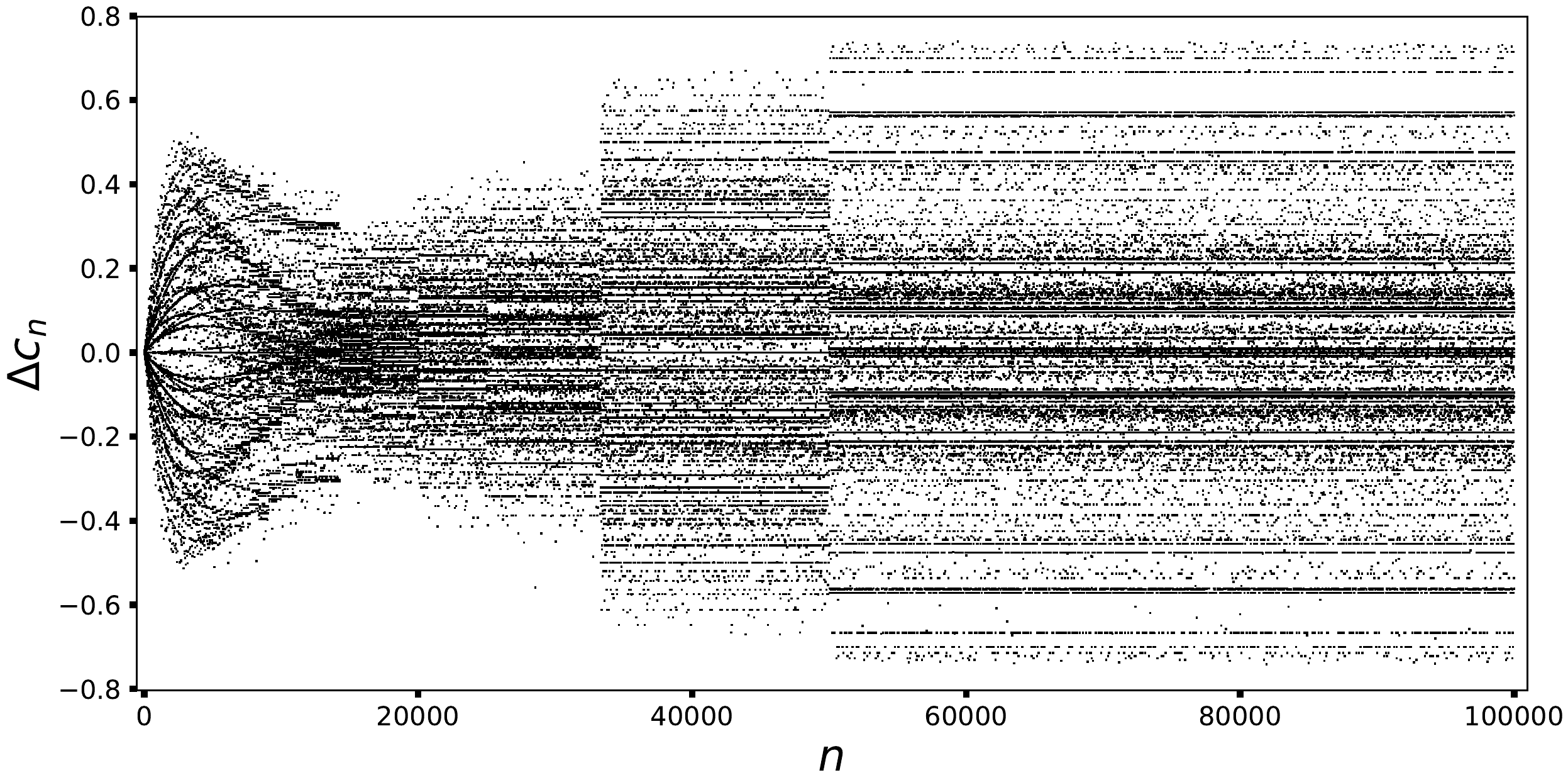}
\label{fig:10a}
\end{subfigure}
\vspace*{\fill} 
\begin{subfigure}{1\textwidth}
\includegraphics[width=1\textwidth]{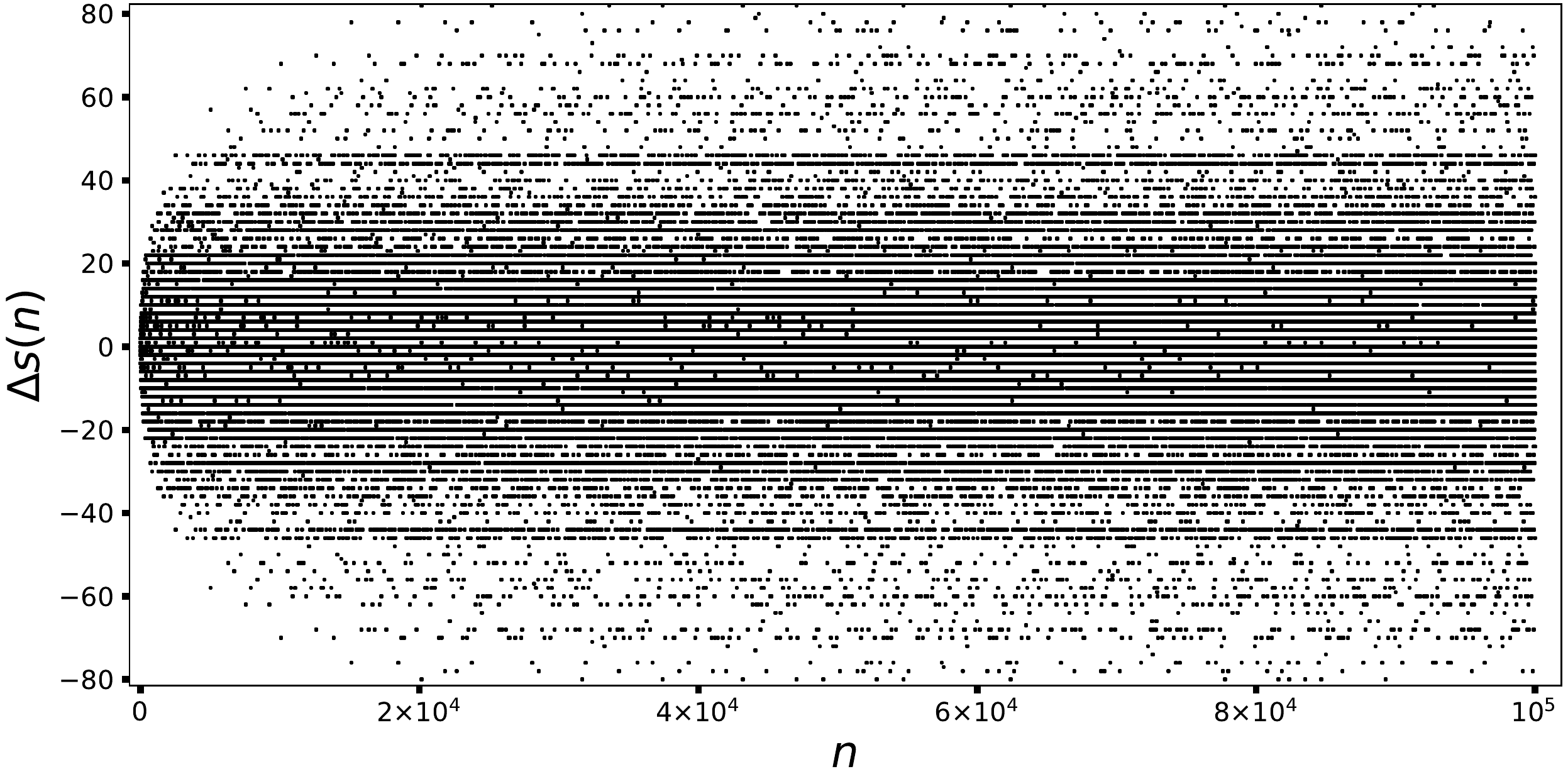}
\label{fig:10b}
\end{subfigure}
\caption{(a) Plot of $ \Delta c_n = c_n - c_{n+1}$  vs nodes $n$ for $N=10^5$. We observe that there are parallel lines and characteristic jumps along with a statistical symmetry about $\Delta c_n = 0$. (b) Plot of $ \Delta s(n) = s(n) - s(n+1)$  vs nodes $n$ for $N=10^5$.  There are also parallel lines and symmetry about $\Delta s(n) = 0$ as in (a) but no jumps.}\label{fig:10}
\end{figure}\\

\newpage

\bibliographystyle{comnet}
\bibliography{Ambika}

\end{document}